\documentclass{sig-alt-full}

\begin{document}

\conferenceinfo{ISSAC'03,} {August 3--6, 2003, Philadelphia, Pennsylvania, USA.} 
      \CopyrightYear{2003} 
      \crdata{1-58113-641-2/03/0008}

\title{Determining the Automorphism Group of a Hyperelliptic 
Curve
}
%

\numberofauthors{1}

\author{
\alignauthor   Tanush Shaska\titlenote{
Useful suggestions of the anonymous referee are gratefully acknowledged.} \\
    \affaddr{Department of Mathematics}\\   
       \affaddr{University of California at Irvine}\\
    \affaddr{103 MSTB, Irvine, CA 92697}\\ 
}

\date{}

\newtheorem {Theorem}{Theorem}[section]
\newtheorem {prop}[Theorem]{Proposition}
\newtheorem {definition}[Theorem]{Definition}
\newtheorem {Example}[Theorem]{Example}
\newtheorem {Exercise} [Theorem]{Exercise}
\newtheorem {remark}[Theorem]{Remark}
\newtheorem {Corollary}[Theorem]{Corollary}
\newtheorem {lemma}[Theorem]{Lemma}
\newtheorem {Result}{Result}

\maketitle

\def\Z{\mathbb Z}
\def\Q{\mathbb Q} \def\C{\mathbb C} \def\bP{\mathbb P}
\def\cB{\mathcal B} \def\cA{\mathcal A} \def\L{\mathcal L}
\def\cR{\mathcal R} \def\H{\mathcal H} \def\M{\mathcal M}
\def\N{\mathcal N} \def\w{\widetilde} \def\l{\lambda} \def\s{\sigma}

\def\Qu{\mathcal Q}

\def\a{\alpha} \def\aa{\bar \alpha} \def\b{\beta} \def\p{\mathfrak p}
\def\P{\mathcal P} \def\e{\varepsilon} \def\iso{\equiv}
\def\sem{{\rtimes}}

\def\o{\oplus}

\def\g{\gamma} \def\bg{\bar \gamma} \def\u{\mathfrak u}
\def\U{\mathfrak U} \def\k{\bar k} \def\iso{{\, \cong\, }}
\def\nor{{\,  \vartriangleleft \, }} \def\<{\langle} \def\>{\rangle}
\def\emb{\hookrightarrow } \def\Aut{\mbox {Aut}}
\def\bAut{\overline{\mbox {Aut}}} \def\D{\mathcal D}

\def\G{\bar G}

\def\ff{\bar \psi} \def\f{\psi} \def\ff{\phi} \def\f{\psi}
\def\X{\mathcal X} \def\Y{\mathcal Y} \def\fH{\mathfrak H} \def\t{\mu}
\def\d{\delta} \def\m{\mu} \def\v{\mathfrak v} \def\r{\mu}
\def\bk{\overline k}

\begin{abstract}
In this note we discuss techniques for determining the automorphism
group of a genus $g$ hyperelliptic curve $\X_g$ defined over an
algebraically closed  field $k$ of characteristic zero.  The first
technique  uses the classical $GL_2 (k)$-invariants of binary
forms. This is a practical method for curves of small genus, but has
limitations as the genus  increases, due to the fact that  such
invariants are not known for large genus.

The second approach, which uses dihedral invariants of hyperelliptic
curves,  is a very convenient method and works well in all genera.
First we define the normal decomposition of a hyperelliptic curve with
extra  automorphisms.  Then dihedral invariants are defined in terms
of the coefficients of this normal decomposition. We define such
invariants independently of the automorphism group $\Aut (\X_g)$.
However, to compute such invariants the curve is required to be in its
normal form. This requires solving a nonlinear system of equations.

We find conditions in terms of classical invariants of binary forms for a curve 
to have reduced automorphism group $A_4$, $S_4$, $A_5$. As far as we are aware,
such results have not appeared before in the literature.

\end{abstract}





\medskip

\section{Introduction}

Let $\X_g$ be an algebraic curve of genus $g$ defined over an
algebraically closed  field $k$ of characteristic zero.   We denote by
$\Aut(\X_g)$  the group of analytic (equivalently, algebraic)
automorphisms of $\X_g$. Then $\Aut(\X_g)$ acts on the finite set of
Weierstrass points of $\X_g$.  This action is faithful unless $\X_g$
is hyperelliptic, in which case its kernel is the group of order 2
containing the hyperelliptic involution of $\X_g$. Thus in any case,
$\Aut (\X_g)$ is a finite group. This was first proved by Schwartz.
In 1893   Hurwitz discovered what is now called the Riemann-Hurwitz
formula.  From this he derived that
$$| \Aut (\X_g)| \ \ \leq \ \ 84 \, \, (g-1)$$ which is  known as the
Hurwitz bound.  However, it is not an easy task to compute the
automorphism group of a given  algebraic curve. Even compiling a list
of possible candidates for a small genus $g$  is quite difficult. In
\cite{MS} we provide an algorithm which computes such lists.  We give
a complete list for $g=3$ and list ``large'' groups for $g\leq
10$. This work is based on previous work of Breuer, among many others;
see \cite{MS} for a complete list of references.

If $\X_g$ is hyperelliptic then $\Aut(\X_g)$ is a degree 2 central
extension of $\Z_n, D_n, A_4, S_4, A_5$. We will explain this  briefly
in section 2. However, computing $\Aut(\X_g)$ for a given $\X_g$ is
still difficult. Even sophisticated computer algebra packages do not
have such capabilities for $g \geq 3$.  The case $g=2$ has recently
been implemented in Magma \cite{Magma} and is based  on methods used
in \cite{SV1}.

In this short note we will focus on determining $\Aut(\X_g)$ for a
given  genus $g$ hyperelliptic curve $\X_g$. We will not prove any of
the results.  The interested reader can check  \cite{GS}, 
\cite{Sh3}, \cite{Sh5},  \cite{Sh6}, \cite{Sh7}, or  \cite{SV1} for
details. Most of the papers above have focused on studying the locus of 
all hyperelliptic curves of  genus $g$ whose automorphism group contains 
a subgroup  $G$ and  inclusions between such loci. 
In this paper we combine the above results to get a treatment for
all hyperelliptic curves in all genera. 
We generalize the notion  of  dihedral invariants of hyperelliptic
curves with extra involutions discovered in \cite{GS} for all hyperelliptic
curves with extra automorphisms (cf. Theorem 5.1.). Using these dihedral 
invariants  and  classical
invariants of binary forms of degree $2g+2$ we  discover some nice
necessary conditions for  a curve to have reduced automorphism group
$A_4, S_4, A_5$ (cf. section 5).  


\smallskip

\noindent
{\bf Notation:} We will use the term ``curve'' to mean a ``compact
Riemann surface''. Throughout this paper $\X_g$ denotes a
hyperelliptic curve of genus $g\ge 2$.  
 $D_n$ denotes the dihedral group of order $2n$.


\section{Preliminaries}
Let $k$  be an algebraically  closed field of characteristic  zero and
$\X_g$  be a  genus  $g$  hyperelliptic curve  given  by  the equation
$Y^2=F(X)$, where $\deg(F)=2g+2$. Denote  the function field of $\X_g$
by $K:=k(X,Y)$.  Then, $k(X)$ is the  unique degree 2 genus zero
subfield of  $K$. $K$ is  a quadratic extension field of $k(X)$
ramified exactly at $d=2g+2$ places $\a_1, \dots , \a_d$  of $k(X)$.
The corresponding places  of $K$ are called the {\it Weierstrass
points} of $K$.   Let $\P:=\{ \a_1, \dots , \a_d \}$ and $G=Aut(K/k)$.
Since $k(X)$  is the only  genus 0 subfield  of degree  2  of $K$,
then  $G$  fixes $k(X)$.  Thus, $G_0:=Gal(K/k(X))=\< z_0 \>$, with
$z_0^2=1$, is central in $G$. We call  the {\it reduced automorphism
group}  of $K$ the  group $\G:=G/G_0$.  By a theorem of Dickson, $\G$
is isomorphic to one of the following: 
$$ \Z_n, D_n, A_4, S_4, A_5,$$ 
with branching indices of the corresponding cover $\ff :\bP^1
\to \bP^1/ \G$ given respectively by
\begin{equation}\label{red_sig}
(n,n), (2, 2, n), (2, 3, 3), (2, 4, 4), (2, 3, 5).
\end{equation}

In \cite{BS}  all subgroups of $G$ are classified and in \cite{Bu} all groups that occur as full automorphism groups of hyperelliptic curves are classified.  We use  the notation of \cite{Bu} and define $V_n$, $H_n$, $G_n$, $U_n$,  $W_2$, $W_3$   as follows:

\begin{equation}
\begin{split}
V_n := & \<  \, \, x, y \, | \, \, x^4, y^n,  (x y)^2, (x^{-1}y)^2 \>, \\  
H_n := & \<\,  x, y \, \, | \, \, x^4, y^2x^2, (xy)^n \, \>,   \\ 
G_n:= & \< \, x, y \, \, | \, \, x^2 y^n, y^{2n},  x^{-1} y x y \, \>,\\
U_n := & \< \, x, y \, | \, \, x^2, y^{2n},  xyxy^{n+1} \>,  \\ 
W_2:=  & \< \, x, y \, | \, \, x^4, y^3, y x^2  y^{-1} x^2, (x y)^4 \>, \\ 
W_3:=  & \< \, x, y \, | \, \, x^4, y^3, x^2 (x y)^4, (x y)^8 \> \\
\end{split}
\end{equation}

The following is proven in \cite{Bu}.

\begin{Theorem}
The automorphism group of a hyperelliptic curve is one of the  following $D_n$, $\Z_n$, $V_n$, $H_n$, $G_n$, $U_n$, $GL_2 (3)$, $W_2$, $W_3$.
\end{Theorem}

The reader should be careful when reading Theorem 3.1., in  \cite{Bu}. It seems
as the cases $H_1$ and $G_1$ (which are isomorphic to $ \Z_4$) must be excluded. 
For example, for $g=2$,  according to Theorem 3.1., $H_1\iso \Z_4$ must occur
as an automorphism group, but it is well known that this is not the case; see \cite{SV1} 
among many others. It is safe to exclude these cases since the group is cyclic and 
corresponds to case 3 of Table 1. 

Also, for  $g=3$ let   $N=3$ in the case 3.d, of Table 2 in \cite{Bu}. 
This case is not excluded  from Theorem 3.1., (pg. 273). In this case the group is
$D_3$ (dihedral group of order 6) and this group does not occur as an automorphism group
of a genus 3 hyperelliptic curve; see \cite{MS}.


\subsection{MODULI SPACES  OF COVERS}

\def\bC{{\bf C}}

Let $\phi_0: \X_g \to \bP^1$ be the cover which corresponds to the
degree 2  extension $K/k(X)$.  Then, $\psi:= \phi \circ \phi_0$ has
monodromy group $G:=\Aut(\X_g)$.  From basic covering theory, the
group  $G$ is embedded in the group $S_n$, where $n=\deg \psi$. There
is an $r$-tuple $\bar \s := (\s_1, \dots , \s_r)$, where $\s_i\in S_n$
such that  $\s_1, \dots , \s_r$ generate $G$ and $\s_1 \cdots \s_r=1$.
The signature of $\psi$ is an $r$-tuple of conjugacy classes  $\bC
:=(C_1, \dots , C_r)$ in $S_n$ such that $C_i$ is the conjugacy class
of $\s_1$.  We use the notation $n^p$ to denote the conjugacy class of
permutations which are a product of $p$ cycles of length $n$.  Using
the signature of $\phi: \bP^1 \to \bP^1$ given in \eqref{red_sig} and
the Riemann-Hurwitz formula, one  finds  out the signature  of $\psi:
\X_g \to \bP^1$ for any given  $g$ and $G$. A natural question is if a given 
group $G$ occurs as an automorphism group of a curve $\X_g$ with more then
one signature $\bC$ (cf. Theorem \ref{thm_2_2}).

For a fixed $G, \bC$ the family of covers $\psi : \X_g \to \bP^1$ is a
Hurwitz space $\H (G, \bC)$. This is a quasiprojective variety, not a
priori connected.  To show irreducibility of $\H (G, \bC)$ one has to
show that there is only one braid  orbit in the set of Nielsen classes
$Ni (G, \bC)$.

There is a map $$\Phi : \H (G, \bC) \to \H_g$$ where $\H_g$ is the
moduli space  of genus $g$ hyperelliptic curves. We denote by $\d (G,
\bC)$ the dimension in $\H_g$ of $\Phi (\H (G, \bC))$.  Further $i(G)$
denotes the number of involutions of $G$.

\begin{Theorem}\label{thm_2_2}
For each  $g\geq 2$, the groups $G$  that occur as automorphism groups 
 and their signatures $\bC$  are given in Table 1. Moreover;
$\H (G, \bC)$ is an irreducible algebraic variety of dimension  $\d
(G, \bC)$ as given in Table 1.
\end{Theorem}

\begin{tiny}
\begin{table*}[hbt!] \label{table}
      \vskip 0.2cm 
\begin{center} 
\renewcommand{\arraystretch}{1.24}
\begin{tabular}{||c|c|c|c|c|c|c||} 
\hline \hline  
$G$ & $\G$ & $\d (G, \bC)$ & $\d, \, n, \,g  $ & $\bC=(C_1, \dots C_r)$
& $ \ff: \bP^1\to \bP^1$  &i(G)  \\  
\hline \hline  $\Z_2 \o
\Z_n $ & &$\frac {2g+2} n -1$&$n < g+1$ & $(n^2, n^2, 2^n, \dots ,
2^n)$ &&   3 \\  $\Z_{2n}$       &$\Z_n$  &$\frac {2g+1}n -1$  &  &
$(n^2, 2n, 2^n, \dots , 2^n )$& $(n, n)$& 1\\  $\Z_{2n}$       &
&$\frac {2g} n -1 $ &$n < g$ & $(2n, 2n, 2^n, \dots , 2^n)$  & &   1
\\  \hline \hline  $\Z_2 \o D_n$   &        &$\frac {g+1} n$ &  & $(
n^4, 2^{2n}, \dots , 2^{2n} )$ & &  2n+3\\ $V_n$  &  &$\frac
{g+1}n-\frac 1 2$    &  & $(n^4, 4^n, 2^{2n}, \dots , 2^{2n} )$&  &
n+3\\  $D_{2n}$ &$D_n$&$\frac g n$            &  & $((2n)^2, 2^{2n},
\dots , 2^{2n} )$& $(2^n, 2^n, n^2)$ &n+1  \\  $H_n$ &     &$\frac
{g+1} n -1$     &$n < g+1$  &$(4^n, 4^{n}, n^4, 2^{2n} \dots , 2^{2n}
)$ &  &  3 \\  $U_n$    &     &$\frac g n- \frac 1 2$ &$g \neq 2$
&$(4^n, (2n)^2, 2^{2n}, \dots , 2^{2n})$ &   & n+1 \\ $G_n$    &
&$\frac g n - 1$        &$n < g$  &$(4^n, 4^n, (2n)^2, 2^{2n}, \dots ,
2^{2n})$    &   &   1   \\ 
\hline \hline 
$\Z_2\o A_4$ &  &$\frac {g+1}
6$&  &$( 3^8, 3^8, 2^{12}, \dots , 2^{12} )$ & &  \\ 
$\Z_2\o A_4$ &
&$\frac {g-1} 6$&  & $( 3^8, 6^4, 2^{12}, \dots , 2^{12} )$& & 7 \\
$\Z_2\o A_4$&$A_4$&$\frac {g-3} 6$&$\d\neq 0$ & $( 6^4, 6^4, 2^{12},
\dots , 2^{12})$&$(2^6, 3^4, 3^4)$ & \\
$SL_2(3) $ &  & $\frac {g-2} 6$ &$\d\neq 0$ &$(4^6, 3^8, 3^8, 2^{12},
\dots , 2^{12})$ &  & \\ 
$SL_2(3) $ & & $\frac {g-4} 6$ & &$(4^6, 3^8,
6^4, 2^{12}, \dots , 2^{12})$ & &  1 \\ 
$SL_2(3) $ &  & $\frac {g-6}
6$ &$\d\neq 0$ &$(4^6, 6^4, 6^4, 2^{12}, \dots , 2^{12})$ & &  \\
\hline \hline $\Z_2\o S_4$ &  &$\frac {g+1} {12}$& &$( 3^{16}, 4^{12},
2^{24}, \dots , 2^{24} )$& & \\ 
$\Z_2\o S_4$ &  &$\frac {g-3} {12}$&
&$( 6^{8}, 4^{12}, 2^{24}, \dots , 2^{24} )$& & 19\\ 
$GL_2(3)$ &  &$\frac
{g-2} {12}$&    &$( 3^{16}, 8^{6}, 2^{24}, \dots , 2^{24} )$& & \\
$GL_2(3)$&$S_4$&$\frac {g-6}{12}$&       &$( 6^{8}, 8^{6}, 2^{24}, \dots ,
2^{24} )$&$(2^{12}, 3^8, 4^6)$&13 \\ 
$W_2$ &  &$\frac {g-5} {12}$& &$(
4^{12}, 4^{12}, 3^{16},  2^{24}, \dots , 2^{24} )$& & 7\\ $W_2$ &
&$\frac {g-9} {12}$& &$( 4^{12}, 4^{12}, 6^{8}, 2^{24}, \dots , 2^{24}
)$& & \\ $W_3$ & &$\frac {g-8} {12}$ & &$( 4^{12}, 3^{16}, 8^{6},
2^{24}, \dots , 2^{24} )$& & 1 \\ $W_3$ & &$\frac {g-12} {12}$& &$(
4^{12}, 6^{8}, 8^{6}, 2^{24}, \dots , 2^{24} )$& & \\ \hline \hline
$\Z_2\o A_5$ & & $\frac {g+1} {30}$& &$(3^{40}, 5^{24}, 2^{60}, \dots
, 2^{60})$&    &\\ $\Z_2\o A_5$ & & $\frac {g-5} {30}$& &$(3^{40},
10^{12}, 2^{60}, \dots , 2^{60})$&     & 31 \\ $\Z_2\o A_5$ & & $\frac
{g-15} {30}$& &$(6^{20}, 10^{12}, 2^{60}, \dots , 2^{60})$&   &\\
$\Z_2\o A_5$ & $A_5$ & $\frac {g-9} {30}$& &$(6^{20}, 5^{24}, 2^{60},
\dots , 2^{60})$& $( 2^{30}, 3^{20}, 5^{12} )$   &\\ $SL_2(5)$&
&$\frac {g-14} {30}$& &$(4^{30}, 3^{40}, 5^{24}, 2^{60}, \dots ,
2^{60})$&   &\\ $SL_2(5)$&  &$\frac {g-20} {30}$& &$(4^{30}, 3^{40},
10^{12}, 2^{60}, \dots , 2^{60})$&    & 1\\ $SL_2(5)$&  &$\frac {g-24}
{30}$ & &$(4^{30}, 6^{20}, 5^{24}, 2^{60}, \dots , 2^{60})$&   &\\
$SL_2(5)$&  &$\frac {g-30} {30}$ & &$(4^{30}, 6^{20}, 10^{12}, 2^{60},
\dots , 2^{60})$&   &\\ \hline \hline
\end{tabular}
\end{center}
\caption{Automorphism groups of hyperelliptic curves}
\end{table*}

\end{tiny}


Finding algebraic descriptions of Hurwitz spaces is in general a difficult
problem. In \cite{Sh7} it is shown that each of the spaces $\H (G, \bC)$ is a 
rational variety. Further, the inclusions between such loci are studied. 

Let $t$ be the order of an automorphism of an algebraic curve $\X_g$
(not necessary hyperelliptic). Hurwitz \cite{Hu} showed that $t \leq
10 (g-1)$.  In 1895, Wiman improved this bound to be $t \leq 2 (2g+1)$
and showed that it is the best possible. Thus, if a cyclic group $H$
occurs as an automorphism group then  $|H| \leq 2 (2g+1)$. Indeed,
this bound can be achieved for any genus via a hyperelliptic
curve. For example, the curve $$Y^2=X(X^{2g+1}-1)$$ has automorphism
group  the cyclic group of order $4g+2$.  This is the second case in
Table 1, when $n=2g+1$. The family of such curves is 0-dimensional in
$\H_g$.

Now we turn our attention to determining if  a given  curve $\X_g$
belongs to any of the families of Table 1. In other words, find
conditions in terms of the  coefficients of $\X_g$ such that $\X_g$
belong to a family in Table 1.  This would determine the $\Aut
(\X_g)$.

\section{Invariants of Binary Forms}

In this section we define the action of $ GL_2(k)$ on binary forms and
discuss the basic notions of their invariants.  Let $k[X,Z]$  be the
polynomial ring in  two variables and  let $V_d$ denote  the
$(d+1)$-dimensional  subspace  of  $k[X,Z]$  consisting of homogeneous
polynomials.
\begin{equation}  \label{eq1}
f(X,Z) = a_0X^d + a_1X^{d-1}Z + ... + a_dZ^d
\end{equation}
of  degree $d$. Elements  in $V_d$  are called  {\it binary  forms} of
degree $d$.  We let $GL_2(k)$ act as a group of automorphisms on $
k[X, Z] $   as follows:
\begin{equation}
 M =
\begin{pmatrix} a &b \\  c & d
\end{pmatrix}
\in GL_2(k),   \textit{   then       }
\quad  M  \begin{pmatrix} X\\ Z \end{pmatrix} =
\begin{pmatrix} aX+bZ\\ cX+dZ \end{pmatrix}.
\end{equation}
This action of $GL_2(k)$  leaves $V_d$ invariant and acts irreducibly
on $V_d$.
%
%
Let $A_0$, $A_1$,  ... , $A_d$ be coordinate  functions on $V_d$. Then
the coordinate  ring of $V_d$ can be  identified with $ k[A_0  , ... ,
A_d] $. For $I \in k[A_0, ... , A_d]$ and $M \in GL_2(k)$, define $I^M
\in k[A_0, ... ,A_d]$ as follows
\begin{equation} \label{eq_I}
{I^M}(f):= I(M(f))
\end{equation}
for all $f \in V_d$. Then  $I^{MN} = (I^{M})^{N}$ and Eq.~(\ref{eq_I})
defines an action of $GL_2(k)$ on $k[A_0, ... ,A_d]$.
A homogeneous polynomial $I\in k[A_0, \dots , A_d, X, Z]$ is called a
{\it covariant}  of index $s$ if
$$I^M(f)=\delta^s I(f)$$
where $\delta =\det(M)$.  The homogeneous degree in $A_1, \dots , A_n$
is called the {\it degree} of $I$,  and the homogeneous degree in $X,
Z$ is called the {\it  order} of $I$.  A covariant of order zero is
called {\it invariant}.  An invariant is a $SL_2(k)$-invariant on
$V_d$.

We will use the symbolic method of classical theory to construct
covariants of binary forms.    Let
\begin{equation}
\begin{split}
f(X,Z):= & \sum_{i=0}^n
\begin{pmatrix} n \\ i \end{pmatrix}a_i X^{n-i} \, Z^i, \\
 g(X,Z) := & \sum_{i=0}^m   \begin{pmatrix} m \\ i \end{pmatrix} b_i
 X^{n-i} \, Z^i    \\
\end{split}
\end{equation}
be binary forms of  degree $n$ and $m$ respectively with coefficients
in $k$. We define the {\bf r-transvection}
\begin{equation}
(f,g)^r:= c_k \cdot \sum_{k=0}^r (-1)^k
\begin{pmatrix} r \\ k
\end{pmatrix} \cdot
\frac {\partial^r f} {\partial X^{r-k} \, \,  \partial Y^k} \cdot
\frac {\partial^r g} {\partial X^k  \, \, \partial Y^{r-k}}
\end{equation}
%
where $c_k=\frac {(m-r)! \, (n-r)!} {n! \, m!}$.  It is a homogeneous
 polynomial in $k[X, Z]$ and therefore a covariant of order $m+n-2r$
 and degree 2. In general, the $r$-transvection of two covariants of
 order $m, n$ (resp., degree $p, q$) is a covariant of order $m+n-2r$
 (resp., degree $p+q$).

For the rest of this paper $F(X,Z)$ denotes a binary form of order
$d:=2g+2$ as below
\begin{equation}
F(X,Z) =   \sum_{i=0}^d  a_i X^i Z^{d-i} = \sum_{i=0}^d
\begin{pmatrix} n \\ i
\end{pmatrix}    b_i X^i Z^{n-i}
\end{equation}
where $b_i=\frac {(n-i)! \, \, i!} {n!} \cdot a_i$,  for $i=0, \dots ,
d$.  We denote invariants (resp., covariants) of binary forms by $I_s$
(resp., $J_s$) where the subscript $s$ denotes the degree (resp., the
order).  We define the following covariants and invariants:
\begin{equation}\label{covar}
\begin{split}
  I_2 & :=(F,F)^d, \\  J_{4j}   & :=   (F,F)^{d-2j}, \, \,  j=1, \dots
, g, \\   I_4 & :=(J_4, J_4)^4, \\  I_4'    &:=     (J_8, J_8)^8, \\
I_6 & :=((F, J_4)^4, (F, J_4)^4)^{d-4}, \\ I_6^\prime  &:=((F, J_8)^8,
(F, J_8)^8)^{d-8}, \\  I_6^{``} & :=((F, J_{12})^{12}, (F,
J_{12})^{12})^{d-12}, \\ I_3      &:=(F, J_d)^d, \\ M & :=((F, J_4)^4,
(F, J_8)^8)^{d-10}, \\ I_{12}  & :=(M, M)^8\\
\end{split}
\end{equation}

$GL_2(k)$-invariants are called 
{\it absolute invariants}.  We  define
the following absolute invariants:
\begin{equation}
\begin{split}
& i_1:=\frac {I_4'} {I_2^2},\, i_2:=\frac {I_3^2} {I_2^3},\,
i_3:=\frac {I_6^{``} } {I_2^3}, \,  j_1 := \frac {I_6^{'}} {I_3^2}, \\
&  j_2:= \frac {I_6} {I_3^2},  s_1:=\frac {I_6^2} {I_{12}}, \,
s_2:=\frac {(I_6^{'})^2} {I_{12}}, \,  \v_1:= \frac {I_6} {I_6^{``}},
\\ &   \v_2:=\frac {(I_4^{'})^3} {I_3^4}, \, \,  \v_3:= \frac {I_6}
{I_6^{'}}, \,   \v_4:=\frac {(I_6^{``})^2} {I_4^3}, \,   \v_5:=\frac
{I_6^{``}}{I_6^{'}}
\end{split}
\end{equation}
For a given curve $\X_g$ we denote by $I(\X_g)$ or $i(\X_g)$ the
corresponding  invariants. Two isomorphic hyperelliptic curves  have
the same absolute invariants.
\begin{remark}
It is an open problem to determine the field of invariants of binary
form of degree $d \geq 7$.
\end{remark}

\section{Equations of curves}
In this section we state the equations of curves in each case of Table
1.  For a more detailed treatment of these spaces, including proofs,
the reader can check results in   \cite{Sh6}, \cite{Sh7}. The reader
can also check \cite{BG} where equations for each family  are
computed; however the main goal of the book is to study hyperelliptic
Riemann surfaces with real structures.  In this section  $G$ denotes a
group as in the first column of Table 1, and  $\L_g^G$ the locus of
hyperelliptic genus $g$ curves $\X_g$ such that $G$ is embedded  in
$\Aut(\X_g)$.

\subsection{$\bAut(\X_g)$ is isomorphic to $\Z_n$}
If $\bAut(\X_g)\iso \Z_n$ then $\X_g$ belongs to cases 1, 2, 3 in
Table 1.  These loci were studied in detail in \cite{Sh6}.  The family
of  curves are given below:
\begin{equation}
\begin{split}
Y^2= & X^{nt}+ \dots + a_i X^{n(t-i)} + \dots a_{t-1} X^n +1, \\ Y^2=
& X^{nt}+ \dots + a_i X^{n(t-i)} + \dots a_{t-1} X^n +1,  \\ Y^2= &
X\, (X^{nt} + \dots + a_i X^{n(t-i)} + \dots a_{t-1} X^n +1) \\
\end{split}
\end{equation}
where $t$ is respectively $ \frac {2g+2} n,  \frac {2g+1} n,  \frac
{2g} n$.  To classify these curves (up to isomorphism) we need to find
invariants of the $GL_2(k)$-action on $k(a_1, \dots , a_{t-1})$.  The
following
\begin{equation}
u_i:=  a_1^{t-i} \, a_i \, + \, a_{\d}^{t-i} \, a_{t-i}, \quad for
\quad 1 \leq i \leq \d\\
\end{equation}
are called {\it dihedral invariants} for the genus  $g$ and the tuple
$$\u:=(u_1, \dots , u_\d)$$ is called the {\it tuple of dihedral
invariants}.  It can be checked that  $\u=0$ if and only if
$a_1=a_\d=0$. In this case replacing $a_1, a_\d$ by $a_2, a_{\d-1}$ in
the formula above would give new invariants.  The next theorem shows
that the dihedral invariants generate $k(\L_g^G)$.
\begin{Theorem}
$ \L_g^G$  is a $\d$-dimensional rational variety. Moreover,  $
k(\L_g^G) =k(u_1, \dots , u_\d)$.
\end{Theorem}
If $n=2$ then $G$ is the Klein 4-group. Then  $\L_g^G=\L_g$ where
$\L_g$ is the locus of hyperelliptic curves with extra involutions,
see \cite{GS}.  A nice necessary and sufficient condition is found in
\cite{Sh5} in terms of the dihedral invariants for a curve to have
more than three involutions. More precisely, for such curves the
relation  holds: $$2^{g-1}u_1^2 - u_g^{g+1}=0.$$ 
%
\subsection{$\bAut(\X_g)$ is isomorphic to $D_n$}
%
The dihedral group is generated by
$$D_n=\< \s, \tau | \, \, \s^n=\tau^2=1\>$$ such that
$$  \s(X)=\e_n\, X, \quad \tau(X)=\frac 1 X.$$
Then $\s$  fixes $X= 0, \infty$ and $\tau $ fixes $X=\pm 1$ and
permutes 0 and $\infty$.  We let $$G(X):= \prod_{i=1}^t (X^{2n} + \l_i
X^n +1)$$ Then,
\begin{equation}
\begin{split}
G(X)= & X^{2nt}+ a_1 X^{2nt-n} + \dots + a_t X^{nt} + \\ &  a_{t-1}
      X^{(n-1)t} + \dots  + a_1 X^n +1
\end{split}
\end{equation}
where $a_i, i=1,\dots t$ are polynomials in terms of the symmetric
polynomials  $s_1, \dots , s_t$ of $\l_i$ (i.e., $a_1=s_1, a_2=t+s_2,
a_3=(t-1)s_1+s_3,  a_4:=\left( \overset{t} {n/2}  \right) + (t-2)s_2
+s_4, $  etc.).

Depending on whether  $0, \pm 1, $ and  $\infty$ are Weierstrass
points we get the  equations $Y^2=F(X)$ where
\begin{equation}
\begin{split}
 F(X) = &\,  G(X),  \quad (X^n-1) \cdot G(X), \\ &\,  X \cdot  G(X),
        \quad   (X^{2n}-1) \cdot  G(X)\,\\ &\,     X\, (X^n-1) \cdot
        G(X), \quad  X\, (X^{2n}-1) \cdot G(X) \\
\end{split}
\end{equation}
where $n$ is respectively as in cases 4-9 of Table 1.
\begin{remark}
Notice that in all cases $n$  is even; see Theorem 2.1., in \cite{Bu}.
\end{remark}
The case $Y^2=G(X)$ corresponds to the group $\Z_2 \o D_n$.  If $n=2$,
then this is a special case of $G\iso \Z_2 \o \Z_n$. Indeed,
$$2^{g-1}u_1^2 - u_g^{g+1}=0$$
as expected; see \cite{Sh5} for details.  If $n> 2$ then
$$\u=(u_1, \dots , u_g)=(0, \dots , 0) $$ where $\u=(u_1, \dots , u_g)
$ is as defined in \cite{GS}.

\subsection{$\bAut(\X_g)$ is isomorphic to $ A_4$ }
This case is treated in detail in \cite{Sh6}. Let
\begin{small}
\begin{equation}  \label{G}
G_{i} (X)=X^{12} - \l_i X^{10} - 33 X^8 + 2 \l_i X^6 - 33 X^4 - \l_i
X^2+1,
\end{equation}
\end{small}
for $\l_1^2+108\neq 0$. Denote by
$$G(X):=\prod_{i=1}^\d G_i(X)$$
Then, each family is parameterized as in Table 2.
\begin{table}[ht]     
\begin{center}
\renewcommand{\arraystretch}{1.24}
\begin{tabular}{||c|c|c||}
\hline \hline $G$  &  $\delta$  & Equation $Y^2= $   \\  \hline \hline
$\Z_2\o A_4$ & $\frac {g+1} 6$ & $G(X)$ \\ $\Z_2\o A_4$ & $\frac {g-1}
6$ & $ (X^4+2i \sqrt{3}X^2 +1)\cdot  G(X)$ \\
$\Z_2\o A_4$ &$\frac {g-3} 6$& $ (X^8+14 X^4+1)  \cdot G(X)$   \\
$SL_2(3) $ &   $\frac {g-2} 6$ &  $X (X^4-1)  \cdot G(X) $ \\
$SL_2(3) $ &  $\frac {g-4} 6$ & $X (X^4-1)(X^4+2i \sqrt{3} X^2 +1)
\cdot G(X)$\\
$SL_2(3) $ &   $\frac {g-6} 6$ &  $X(X^4-1) (X^8+14X^4+1)\cdot G(X)$
\\   \hline \hline
\end{tabular} 
\end{center}
\caption{Hyperelliptic curves with $\bAut(\X_g)= A_4$}
\end{table}
The following lemma gives a necessary condition that a curve has
automorphism  group $\Z_2 \o A_4$ or $SL_2(3)$.
\begin{lemma} Let $\X_g$ be a hyperelliptic curve of genus $g$ with 
$\bAut(\X_g)\iso A_4$. Then,  $I_4=0$. Moreover;

i) if $g=4$ then $I_2=I_4=I_4^{'}=I_6^{'}=0$

ii) if $g=5, 9,  12$ then  $I_4=I_6=0$

iii) if $g=7, 10$ then $I_2=I_4=I_4^{'}=I_6^{``}=0$

iv) if $g=8$ then $I_4=0$.
\end{lemma}

\subsection{$\bAut(\X_g)$ is isomorphic to $ S_4$ }
In this case the reduced automorphism group is generated by
$$\s (X)=- \frac {x-1} {x+1}, \quad \tau(X)=i X.$$
We also denote
\begin{equation}
\begin{split}
G_i(X):= & X^{24} + \l X^{20} + (759-4\l)X^{16} + 2 (3\l +\\ & 1288)
        X^{12} + (759-4\l)X^{8}+ \l X^{4}+1, \\ \\ R(X):= & X^{12} -
        33 X^8 - 33 X^4 +1, \\  S(X):= & X^8 + 14 X^4 + 1, \\  T(X):=
        & X^4-1.
\end{split}
\end{equation}
Let $$G(X):=\prod_{i=1}^\d G_i(X)$$ where $\d$ is as in Table 1.
Then, the equations of the curves in each case are $Y^2=F(X)$ where
$F$ is as  below (we suppress  $X$):
$$F=  G, \,   S G, \,   T G, \,   S  T G, \,  R G, \,  R S  G, \, R  T
G, \,  R S   T G.
$$
Similar conditions in terms of the classical invariants as in the
previous case can be obtained in this case also.
\begin{lemma}Let $\X_g$ be a hyperelliptic curve of genus $g$ with 
$\bAut(\X_g)\iso S_4$. Then,  $I_4=0$.
\end{lemma}

\subsection{$\bAut(\X_g)$ is isomorphic to $ A_5$ }
We briefly state the equations here. We denote by $G_i(X)$, $R(X)$, $S(X)$,
$T(X)$ the following:
\begin{tiny}
\begin{equation}
\begin{split}
G_i (X):= &    (\l_i -1) X^{60} - 36\, (19\l_i +29) X^{55}+6
(26239\l_i - 42079) X^{50}\\ &   - 540 (23199\l_i -19343) X^{45} + 105
(737719\l_i - 953143)  X^{40} \\ &   - 72 (1815127\l_i - 145087)
X^{35} - 4 (8302981\l_i +49913771)  X^{30}\\ &  + 72 (1815127\l_i -
145087) X^{25} + 105 (737719\l_i - 953143)  X^{20}\\ &   + 540
(23199\l_i -19343) X^{15} + 6 ( 26239\l_i - 42079 ) X^{10} \\ &   +
36\, (19\l_i +29) X^{5}+(\l_i-1)\\ \\ R(X):= & X^{30}+522X^{25}
-10005X^{20} -10005X^{15} -522X^5+1\\ \\ S(X):= & X^{20} - 228X^{15} +
494X^{10} +228X^5 +1\\ \\ T(X):= & X^{10}+10X-1.\\
\end{split}
\end{equation}
\end{tiny}
As above we let $$G(X):=\prod_{i=1}^\d G_i(X).$$ In the order of Table
1 equations are given as $Y^2=F(X)$ where $F(X)$ is as given  as (we
suppress  $X$):
$$F=  G, \,   S G, \,   T G, \,   S  T G, \,  R G, \,  R S  G, \, R  T
G, \,  R S   T G
$$

These curves can be expressed as $Y^2 =\, M(X^2)$ or $Y^2=X\cdot
M(X^2)$ where  $M$ is a polynomial in $X^2$. This fact will be used in
the next section.  The expressions are rather large and we will not
state them here.  However, we get the following useful fact:
\begin{lemma}Let $\X_g$ be a hyperelliptic curve of genus $g$ with 
$\bAut(\X_g)\iso A_5$. Then,
$I_{4}=I_4^\prime=I_6=I_6^\prime=I_{12}=0$.
\end{lemma}

\section{Determining the automorphism group of a given curve}
Let $\X_g$ be given. We want to determine $\Aut(\X_g)$.  In order to
find an algorithm which would work for any $g$ we would have to  check
whether $\X_g$ can be written in any of the forms above. Thus, we want
to find if there is a coordinate change
$$X \to \frac {aX + b} {cX+d}$$
which transforms $\X_g$ to one of the forms of section 4. This would
require solving a system of equations for each case and therefore
would not be efficient.

\subsection{Using classical invariants}

For a fixed $g$ we know the dimension $\d$ of the locus $\L_g^G$. We
compute enough  absolute invariants to generate this locus. Thus, we
determine the loci $\L_g^G$  for all $G$ in Table 1 in terms of some
invariants $i_1, \dots i_{\d+1}$. These loci are computed only once
for each $g$.  Then,  for a particular curve we simply compute these
invariants and check if they  generate any of the loci $\L_g^G$.
These spaces were computed in detail in \cite{Sh7} for  $\G=A_4$. We
will illustrate with $g \leq 12$ and $\G \iso  A_4, S_4, A_5$.

We define $\p(\X_g)$ as follows:
\begin{eqnarray*}
\p(\X_g):=(\p_1, \p_2)=  \left\{ \aligned \v_1, \qquad \qquad \qquad
\qquad \textit{  if  } g=4,   \\ (i_1, i_2),  \quad \textit{  if  }
g=5, 9, \textit{  and   }   I_2\neq 0   \\ \v_2,   \quad \quad
\textit{  if  } g=5, 9, \textit{  and   }   I_2=0   \\  (j_1, j_2),
\qquad \textit{  if  } g=7, \, \,  \textit{  and   }   I_3 \neq 0 \\
\v_3,  \qquad \quad \textit{  if  } g=7, \, \,  \textit{  and   } I_3
= 0 \\ (i_1, i_3),  \quad \textit{  if  } g=8, 12,  \textit{  and }
I_2\neq 0  \\ \v_4,  \qquad \textit{  if  } g=8, 12,  \textit{ and   }
I_2 = 0  \\ (s_2, s_1),  \quad \textit{  if  } g=10, \, \, \textit{
and   }   I_{12}\neq 0     \\ \v_5,  \qquad \textit{  if  } g=10, \,
\,  \textit{  and   } I_{12}= 0 \\ \endaligned \right.
\end{eqnarray*}
From Lemma 4.3. and results for cases $\G \iso S_4, A_5$ one can check
that $\p(\X_g)$ is well defined.  Moreover, the subvariety $\L_g^G$ is
1-dimensional if $\G$ is  isomorphic to $A_4, S_4, A_5$.   For each
parametric curve $\X_g$ of the previous section we compute $\p(\X_g)$
in terms of  the parameter $\l$. Eliminating $\l$ gives an equation
for $\L_g^G$, see \cite{Sh7} for explicit equations.

The following algorithm determines if the automorphism group of a
hyperelliptic genus $g\leq 12$ curve is isomorphic to $\Z_2 \o A_4,
\Z_2\o S_4, \Z_2\o A_5, SL_2(3), SL_2(5), GL_2(3), W_2, W_3$.

\medskip

\noindent {\sc Algorithm 1:}

\smallskip

{\bf Input:} A hyperelliptic curve $\X_g: Y^2= F(X, Z)$.

{\bf Output:} Determine if the automorphism group $\Aut(\X_g)$
is one of  $\Z_2 \o A_4$, $\Z_2\o S_4$, $\Z_2\o A_5$, $SL_2(3)$, $SL_2(5)$,
$GL_2(3)$, $W_2$, $W_3$.

\medskip

{\bf Step1:} Compute $I_4 (\X_g)$. If $I_4\neq 0$ then $\Aut(\X_g)$ is
not isomorphic to any of $\Z_2 \o A_4, \Z_2\o S_4,$ $\Z_2\o A_5$,
$SL_2(3), SL_2(5)$, $GL_2(3), W_2, W_3$. Otherwise go to {\bf Step 2}.

{\bf Step 2:} Compute $\p(\X_g)$.

{\bf Step 3:} Find  $\L_g^G$ which
is satisfied by  $\p(\X_g)$ (equations are given in \cite{Sh6}). Then,
$\Aut(\X_g)$ is isomorphic to $ G$.

\medskip

The definition of $\p(\X_g)$ is a little more elaborate for $\G=\Z_n,
D_n$ since the dimension of $\L_g^G$ is $> 1$. Once the definition of
the moduli point is  modified and the corresponding $\L_g^G$ are
computed the following can be used:

\medskip

\noindent {\sc Algorithm 2:}

\smallskip

{\bf Input:} A hyperelliptic curve $\X_g: Y^2= F(X, Z)$.

{\bf Output:} The automorphism group $\Aut(\X_g)$.

\medskip

{\bf Step 1:} Compute $\p(\X_g)$.

 {\bf Step 2:} Find  $\L_g^G$ which
is satisfied by  $\p(\X_g)$. Then, $\Aut(\X_g)$ is isomorphic to $ G$.

\medskip

The above method of classical invariants is difficult to implement for
large $g$.  That's because finding enough absolute invariants is not
an easy task for large $g$.  Also the expressions of these invariants
and the equations for the loci $\L_g^G$  get very large as $g$ grows.
In order to deal with these problems we use the  dihedral invariants
which will be    explained next.

\subsection{Using dihedral invariants}

In section 4.1., we introduced dihedral invariants for hyperelliptic curves
$\X_g$ such that $\Aut (\X_g) \iso \Z_n$. In this section we generalize this
approach to all hyperelliptic curves with extra automorphisms. 
Theorem 5.1., makes this generalization possible.

Let $\X_g$ be an hyperelliptic curve with extra automorphisms. The
following lemma gives a general description of how to write an
equation for $\X_g$.

\noindent
\begin{Theorem}\label{thm5}
 Let $\X_g$ be a hyperelliptic curve  with $$| \Aut(\X_g)| > 2.$$ 
Then, $\X_g$ can be written as
\begin{equation}\label{normal}
  Y^2=F(X^n), \quad \textit{ or } \quad Y^2=X\cdot F(X^n),
\end{equation}
where $n=2$ or $n$ is odd and divides $2g+2, 2g+1, g$.   Moreover, if
$n> 2$ then $\Aut(\X_g)$ is a cyclic group.
\end{Theorem}

Let $\X_g$ be a hyperelliptic curve  with $| \Aut(\X_g)| > 2$ and
written  as in \eqref{normal}. We call this form a {\bf decomposition}
of $\X_g$. Let  $s$ be the smallest $n$ that such decomposition is
possible. Then,
\begin{equation}\label{normal2}
Y^2= F(X^{s}), \quad or \quad Y^2=X \cdot F(X^{s})
\end{equation}
is called the {\bf normal decomposition} or the {\bf normal form} of
$\X_g$ and $s$ is called the {\bf degree } of the decomposition. If no
such decomposition is possible then we say that $s=1$.  Let $\X_g$ be
in its normal decomposition given below:
\begin{equation}
\begin{split}
Y^2= & X^{nt}+ \dots + a_i X^{n(t-i)} + \dots a_{t-1} X^n +1,  \\ Y^2=
& X\, (X^{nt} + \dots + a_i X^{n(t-i)} + \dots a_{t-1} X^n +1) \\
\end{split}
\end{equation}
where $nt=2g+2, 2g+1, 2g$.

We define the following
\begin{equation}
u_i:=  a_1^{t-i} \, a_i \, + \, a_{\d}^{t-i} \, a_{t-i}, \quad for
\quad 1 \leq i \leq \d=t-1,\\
\end{equation}
which are called {\it dihedral invariants} for genus  $g$ and the tuple
$$\U^1:=(u_1, \dots , u_\d)$$ is called the {\it tuple of dihedral
invariants}.  It can be checked that  $\u=0$ if and only if
$a_1=a_\d=0$. Then, let $(a_j, a_{\d-j+1})$ be the first nonzero
tuple.  Replacing $a_1, a_\d$ by $a_j, a_{\d-j+1}$ in the formula
above would give new invariants. Thus, we define
\begin{equation}\label{u_i}
u_i^j:=  a_j^{\d-i+1} \, a_i \, + \, a_{\d-j}^{\d-i+1} \, a_{\d-i+1},
\end{equation}
for $1 \leq i \leq \d$,  and  $ 1 \leq j \leq [\frac {\d+1} 2]$. Then
\begin{equation}
\U^j:=(u_1^j, \dots , u_m^j)
\end{equation}
where $m=\d - 2j$.

\medskip

\noindent {\sc Algorithm 3:}

\smallskip

{\bf Input:} A hyperelliptic curve $\X_g: Y^2= F(X, Z)$.

 {\bf Output:} The automorphism group $\Aut(\X_g)$.

\medskip

{\bf Step 1:} Check whether the curve has a normal decomposition. If
 ``Yes'' then go to {\bf Step 2} otherwise  $\Aut (\X_g)=\Z_2$ 

 {\bf  Step 2:} Compute the degree $s$ of the normal decomposition. If $s$
 is odd then  $\Aut(\X_g) \iso \Z_{2s}$, otherwise go to {\bf Step 3}.

 {\bf Step 3:} Compute the  dihedral invariants $\U_i^j$ of the
 normal decomposition.  Go to {\bf Step 4}.  

 {\bf Step 4:} Find
 $\L_g^G$ which is satisfied by  $\U_i^j$. Then, $\Aut(\X_g)$ is isomorphic 
to $G$.

\medskip

The above method was used in \cite{SV1} and \cite{GS} to determine the
automorphism  group of genus 2 and 3. It has the advantages that it
can be used for any $g$ no matter how large. A disadvantage is that a
nonlinear  system of equations must be solved in order to determine
the normal decomposition.





\begin{Example}
 For genus 2, the curve can be written as
$$Y^2=X^6+a_1 X^4 + a_2 X^2 +1 $$ and its the  dihedral invariants are
$$u_1=a_1^3 +a_2^3, \quad   \quad  u_2 = 2 a_1 a_2, $$
Then,

a) $G\iso V_6$ if and only if  $(u_1, u_2 )=(0,0)$ or 
$$(u_1, u_2 )=(6750, 450).$$

b) $G\iso GL_2(3) $ if and only if $(u_1, u_2) = ( -250, 50)$.

c) $G\iso D_{6}$ if and only if 
$$  u_2^2 - 220 u_2 -16 u_1 +4500=0,$$ 
for $u_2 \neq 18, 140 + 60\sqrt{5}, 50$.

d) $G\iso D_4$ if and only if         $$ 2 u_1^2-u_2^3=0,$$ for
$u_2 \neq 2, 18, 0, 50, 450$. Cases $u_2 = 0, 450, 50 $ are reduced to
cases a),and b) respectively, see \cite{SV1} for details.
\end{Example}
\begin{remark}
The notation used in \cite{SV1} to denote the groups is
different. $V_6$ is this case has order 24 and in \cite{SV1} is
identified as $\Z_3 \sem D_4$.
\end{remark}

\section{Closing remarks}

We briefly described techniques of determining the automorphism group
of a hyperelliptic curve. A combination of both methods sometime
produces  better results. Our goal is to combine these methods and
explicitly compute loci $\L_g^G$ for reasonable $g$ (i.e., $g \leq
60$).

There are polynomial time algorithms to compute the decomposition of a
polynomial $F(X)$ up to an affine transformation $X \to aX +b$, see
\cite{Gu}.  However, this is not sufficient for our purposes since we want
to find such decomposition up to a liner fractional transformation $X
\to \frac {aX + b} {cX + d}$. If a polynomial time algorithm would be
found in this case this would make the second method  preferable
to the first.

Besides computing the automorphism groups the above techniques can
also be used  to answer other questions on hyperelliptic curves.
For example dihedral invariants can be used to determine the field of moduli
of a given curve. The reader can check   \cite{Sh5} for details and open questions
on the field of moduli  and other computational aspects 
of hyperelliptic curves.


\bibliographystyle{abbrv}

%
%

\end{document}